\documentclass{jaums}
\usepackage{amsfonts}
\usepackage{amssymb}
\usepackage{ifthen}
\usepackage{graphicx}
\usepackage{amscd}
\usepackage{amsxtra}
\usepackage{color}

\theoremstyle{thmit} % Numbered and Italic
\newtheorem{thm}{Theorem}[section]
\newtheorem{lem}[thm]{Lemma}
\newtheorem{cor}[thm]{Corollary}

\newtheorem{con}[thm]{Conjecture}
\theoremstyle{thmrm} % Numbered and Roman

\newtheorem*{rem}{Remark}
\newtheorem*{oldproof}{Proof}

\newcounter{alphabet}
\newcounter{tmp}

\def\be{\begin{equation}}
\def\ee{\end{equation}}

\newcommand{\bcon}{\begin{con}}
\newcommand{\econ}{\end{con}}
\newcommand{\A}{{\mathcal A}}

\newcommand{\K}{{\mathcal K}}
\newcommand{\es}{{\mathcal S}}

\newcommand{\R}{{\mathbb R}}

\newcommand{\D}{{\mathbb D}}
\newcommand{\C}{{\mathcal C}}

\newcommand{\N}{{\mathbb N}}

\newcommand{\F}{{\mathcal F}}
\newcommand{\G}{{\mathcal G}}

\newcommand{\beq}{\begin{eqnarray}}
\newcommand{\beqq}{\begin{eqnarray*}}
\newcommand{\eeq}{\end{eqnarray}}
\newcommand{\eeqq}{\end{eqnarray*}}

\newcommand{\ds}{\displaystyle}

\title{Logarithmic coefficients problems in families
related to starlike and convex functions}

\author{Saminathan Ponnusamy
}
\address{S. Ponnusamy\\
 Department of Mathematics\\
Indian Institute of Technology Madras\\
 Chennai-600 036, India.
}
\email{samy@iitm.ac.in}

\author{Navneet Lal Sharma}
\address{N.L. Sharma\\
Department of Mathematics\\
Indian Institute of Technology Madras\\
Chennai-600 036, India.
}

\email{sharma.navneet23@gmail.com}

\author{Karl-Joachim Wirths}
\address{K.-J. Wirths\\
Institut f\"ur Analysis und Algebra\\
TU Braunschweig\\
38106 Braunschweig, Germany.}
\email{kjwirths@tu-bs.de}

\subjclass {Primary: 30C45, 30C20; Secondary: 30C50, 30C55}
\keywords{Univalent function, subordination, Logarithmic coefficients,  Schwarz' lemma
}

\begin{document}

\maketitle
\markboth{S. Ponnusamy, N. L. Sharma, and K.-J. Wirths}{Logarithmic coefficients problems }

\begin{abstract}
Let $\es$ be the family of analytic and univalent functions $f$ in the unit disk $\D$
with the normalization $f(0)=f'(0)-1=0$, and let $\gamma_n(f)=\gamma_n$ denote the logarithmic coefficients of $f\in {\es}$. In this paper, we study bounds for the logarithmic coefficients for certain subfamilies of univalent functions. Also, we consider the families $\F(c)$ and $\G(\delta)$  of functions $f\in {\es}$ defined by
$$ {\rm Re}  \left ( 1+\frac{zf''(z)}{f'(z)}\right )>1-\frac{c}{2}\, \mbox{ and } \,
{\rm Re}  \left ( 1+\frac{zf''(z)}{f'(z)}\right )<1+\frac{\delta}{2},\quad z\in \D
$$
for some $c\in(0,3]$ and $\delta\in (0,1]$, respectively.
We obtain the sharp upper bound for $|\gamma_n|$ when $n=1,2,3$ and $f$ belongs to the classes $\F(c)$ and $\G(\delta)$, respectively. The paper concludes with the following two conjectures:
\begin{itemize}
 \item If $f\in\F (-1/2)$, then  $ \displaystyle |\gamma_n|\le \frac{1}{n}\left(1-\frac{1}{2^{n+1}}\right)$ for $n\ge 1$,
 and
$$ \sum_{n=1}^{\infty}|\gamma_{n}|^{2} \leq \frac{\pi^2}{6}+\frac{1}{4} ~{\rm Li\,}_{2}\left(\frac{1}{4}\right)
   -{\rm Li\,}_{2}\left(\frac{1}{2}\right),
$$
where ${\rm Li}_2(x)$ denotes the dilogarithm function.
  \item If $f\in \G(\delta)$, then  $ \displaystyle |\gamma_n|\,\leq \,\frac{\delta}{2n(n+1)}$ for $n\ge 1$.
\end{itemize}

\end{abstract}

\section{Introduction}

Let $\A$ be the class of analytic functions $f$ defined on the unit disk $\D=\{z\in \mathbb{C}:\,|z|<1\}$ of the form $f(z)=z+ \sum_{n=2}^{\infty}a_nz^n$, and $\es$ denote the subclass of functions $f\in\A$
that are univalent in $\D$. The logarithmic coefficients $\gamma_n$ of $f\in {\es}$ are defined by the formula
\begin{equation}
\label{eq1-sec1}
\log\left (\frac{f(z)}{z} \right )=2\sum_{n=1}^\infty \gamma_n(f)z^n \quad \mbox{ for } z\in \D.
\end{equation}
We use $\gamma_n(f)=\gamma_n$ when there is no confusion. These coefficients play an important role for various estimates in the theory of univalent functions, and some authors use $\gamma_n$ in place of $2\gamma_n$.  Louis de Branges \cite{dB85} showed (see also \cite{AvWir-09} and \cite{FP85}) that for each $n\ge1,$
$$ \sum_{k=1}^n k(n-k+1)|\gamma_n|^2\le \sum_{k=1}^n \frac{n-k+1}k,
$$
where equality holds if and only if $f$ has the form $z/(1-e^{i\theta}z)^2$ for some $\theta\in\R.$
It is known that this proves the famous Bieberbach-Robertson-Milin conjectures about Taylor coefficients of $f\in {\es}$ in its most general form. See \cite{AvWir-09}. For another (shorter) version of de Branges' proof, we refer to \cite{FP87}. Note that for $f(z)=z/(1-e^{i\theta}z)^2$ we have $\gamma_n=e^{ni\theta}/n$ for $n=1,2,\dots.$
The idea of studying the logarithmic coefficients helped Kayumov \cite{Kay05} to solve Brennan's conjecture
for conformal mappings. In this note, we will discuss the logarithmic coefficients problem for  certain subfamilies
of univalent functions and derive some similar inequalities. We now begin with certain preliminaries.

A function $f\in \es$ is called starlike if $f(\D)$ is a domain that is starlike with respect to the origin.
Every starlike function is characterized by the condition ${\rm Re}\,\left(zf'(z)/f(z)\right)>0$ for
$z\in\D$.  A function $f\in \es$ is convex when
the function $g=zf'$ is starlike. A function  $f\in\A$ is said to be close-to-convex
if there exists a real number $\theta$ and a function $g\in \es^*$ such that
${\rm Re}\,\left(e^{i\theta}zf'(z)/g(z)\right)>0$ in $\D$.
We denote by $\es ^*,\, \C$ and $\K$, the class of starlike functions, convex functions and
close-to-convex functions, respectively.
Functions in the class $\K$ of all close-to-convex functions are known to be univalent in $\D$.
The role of $\es$ together with its subfamilies and their importance
concerning strongly starlike functions in geometric function theory
are well documented. See for example~\cite{BK69,Sta66}, the books of Duren \cite{Dur83}
and Goodman \cite{Go}.

If $f\in\es,$ then it is known that $|\gamma_1|\le 1$ and
$$|\gamma_2|=\frac{1}{2} \left|a_3-\frac{1}{2}a_2^2 \right| \le \frac{1}{2}(1+2e^{-2})= 0.635\dots
$$
by using the Fekete-Szeg\"{o} inequality (see~\cite[Theorem 3.8]{Dur83}). For $n \ge 3$,
the logarithmic coefficients problem seems much harder.  The inequality $|\gamma_n| \le 1/n$
holds for $f\in\es^*$ but is not true for the full class $\es$, even in order
of magnitude (see \cite[Theorem 8.4]{Dur83}). Indeed, there exists a bounded function $f\in \es$ with $\gamma_n\ne O(n^{0.83})$.
On the other hand, Roth \cite{Roth07} established the following sharp inequality for $f\in\es$:
\begin{equation}\label{eq:main}
\sum_{n=1}^\infty   p_n|\gamma_n|^2\le \sum_{n=1}^\infty   \frac{p_n}{n^2}, \quad  p_n =\left (\frac{n}{n+1}\right )^2  .
\end{equation}
This inequality is a source for many new inequalities for the logarithmic coefficients of $f\in\es,$ such as
$$\sum_{n=1}^\infty |\gamma_n|^2\le \frac{\pi^2}6.
$$
Elhosh \cite{Elhosh96} proved that if $f\in\K$, then  $|\gamma_n|\le 1/n$.
However, Girela~\cite{Gir2000} pointed out that this bound is false for the class $\K$ when $n\ge 2$.
He proved that there exists a function $f\in\K$ such that $|\gamma_n| > 1/n$ for $n\ge 2$.
It remains an open problem to find the correct order of growth of $|\gamma_n|$ for $f\in\es$ even for $f\in\K$.
In  a recent paper \cite{OPW18}, the authors considered  the inequality of type \eqref{eq:main} and
the order of growth of $|\gamma_n|$ for the  class ${\mathcal G}(c)$, which is defined below, and also
for some other subclasses of $\es$.

Let $\mathcal B$ denote the class  of all analytic functions $\phi$ in $\D$
which satisfy the condition $|\phi(z)|<1$ for $z\in\D$. Functions in $\mathcal B$ are called
Schwarz functions.
Let $f$ and $g$ be two analytic functions in $\D$. We say that $f$ is {\em subordinate} to $g$,
written as $f \prec g$, if there exists a function $\phi\in \mathcal B$ such that $f(z)=g(\phi(z))$
for $ z\in \D.$
In particular, if $g$ is univalent in $\D$, then $f \prec g$ is equivalent to
$f(\D)\subset g(\D)$ and $f(0)=g(0)$. See \cite{Dur83,Pom75}.

\begin{enumerate}
 \item {The class $\es^*(A,B)$} is defined by
$$\es ^*(A,B):=\left \{f\in \A: \, \frac{zf'(z)}{f(z)}\prec \frac{1+Az}{1+Bz} \,
 \mbox{ for } z\in \D \right\},
$$
where $A\in\mathbb{C},-1\le B\le 0$ and $A\neq B$.
The class $\es^*(A,B)$ with the restriction $-1\le B<A\le 1$ has been studied by Janowski~\cite{Jan73}.
In particular, for $B=-1$ and $A=e^{i\alpha}(e^{i\alpha}-2\beta\cos \alpha),\, 0\le\beta <1,$
 the class $\es ^*(A,B)$ reduces to the class of spiral-like functions of order $\beta$, denoted by  $\mathcal{S}_\alpha (\beta )$,
so that
% \item {\bf The Spiral-like class of order $\beta$:}
$$ \mathcal{S}_\alpha (\beta ) %:=\es ^*(e^{i\alpha}(e^{i\alpha}-2\beta\cos \alpha),-1)
= \left \{f\in \A: \, {\rm Re}  \left ( e^{-i\alpha}\frac{zf'(z)}{f(z)}\right )>\beta \cos \alpha\,
 \mbox{ for }  z\in \D\right\},
$$
where $ \beta \in [0,1)$ and $\alpha \in (-\pi/2, \pi/2)$.
Each function in $\mathcal{S}_\alpha (\beta )$ is univalent in $\D$ (see \cite{Lib67}).
Clearly, $\mathcal{S}_\alpha (\beta )\subset \mathcal{S}_\alpha (0)\subset \mathcal{S}$
whenever $0\leq \beta <1$. Functions in $\mathcal{S}_\alpha(0)$  are called \textit{$\alpha$-spirallike}, but they do not necessarily
belong to the starlike family $\mathcal{S}^*$.
The class $\mathcal{S}_\alpha (0)$ was introduced by ${\rm \check{S}}$pa${\rm\check{c}}$ek \cite{Spacek-33}
(see also \cite{Dur83}).

\item {The class $\es \es^*_{\alpha}$ of strongly starlike functions} is defined by (\cite{BK69,Sta66})
$$\es \es^*_{\alpha}:=\left\{f\in \A: \, \frac{zf'(z)}{f(z)}\prec \left(\frac{1+z}{1-z}\right)^{\alpha}
\,  \mbox{ for }  z\in \D \right\},
$$
for $0<\alpha\le 1$.
%This class was introduced by Brannan and Kirwan~\cite{BK69}, and Stankiewicz~\cite{Sta66}.
 For $\alpha=1$, the class $\es \es^*_{\alpha}$
reduces to the class of starlike functions.

%%%%%%%%%%%%%%%%%%%% Class F(C)  %%%%%%%%%%%%%%%%%%%%%%%%%%%
\item {The class $\mathcal{F}(c)$ is defined by}
\begin{align*}
\mathcal{F}(c) & :=\left \{f\in \A: \, {\rm Re}  \left ( 1+\frac{zf''(z)}{f'(z)}\right )>1-\frac{c}{2},
 \, z\in \D \right\}\\
& =\{f\in\A: \, zf'\in \es^*(c-1,-1)\}
\end{align*}
for some $c\in(0,3]$.
If we set $\alpha\,=\,1\,-\,c/2\in [0,1)$, in this choice the family as $\mathcal{F}(c)$ is well-known  and is referred to
as the family of convex functions of order $\alpha$. Clearly, the family $\F(2)$ is the usual
class of normalized convex functions.
In particular, for $c=3$, we have the class $\mathcal{F}(3)$ which attracted the attention of many in the recent years
(see \cite{PonSaYa14} and the references therein).
Also it is important to point out that functions in $\F(3)$ are known to be convex in one direction (and hence, univalent and close-to-convex)
but are not necessarily starlike in $\D$ (\cite{Ume52}).

%%%%%%%%%%%%%%%%  Class G(C)  %%%%%%%%%%%%%%%%%
\item {The class $\G(c)$ is defined by}
\begin{align*}
\G(c)& :=\left \{f\in \A: \, {\rm Re}  \left ( 1+\frac{zf''(z)}{f'(z)}\right )
                   <1+\frac{c}{2}, \, z\in \D \right\}\\
& =\{f\in\A: \, zf'\in \es^*(-(1+c),-1)\}
\end{align*}
for some $c\in (0,1]$.
Set $\G(1)=:\G$. It is known that $\G\subset \es^{*}$ and thus, functions in $\G(c)$ are starlike.
This class has been studied extensively in the recent past, see for instance \cite{OPW13,OPW18}
and the references therein.

\end{enumerate}

%%%%%%%%%%%%%%%%% Main results %%%%%%%%%%%%%%%%%%%%%%%%%%

The paper is organized as follows. In Section \ref{main reults}, we state our main results along with
two conjectures concerning the logarithmic coefficients bound for $\F(c)$ and  $\G(c)$, respectively.
These conjectures have been verified for the first three logarithmic coefficients.
In Section \ref{main reults}, we have proved the
logarithmic coefficients bound completely for the families $\es^*(A,B)$ and $\es \es^*_{\alpha}$ along with inequalities of the type
\eqref{eq:main} for these families. Our final results in Section \ref{main reults} concern the families $\F(c)$ and
$\G(c)$, where we obtain sharp estimates for the initial three coefficients, and non-sharp estimates for fourth and fifth coefficients.
In Section \ref{lemmas}, we recall a few important lemmas which are useful in the sequel.
The proofs of our main results are presented
in Section \ref{proof theorems}.

\section{Main Results}\label{main reults}
%%%%%%%%%%%%%%%%%%%%%%%% theorem-1 %%%%%%%%%%%%%%%%
\begin{thm}\label{thm1}
 For $-1\le B<A\le 1$  and $B\neq 0$, the logarithmic coefficients of $f\in\es^*(A,B)$ satisfy the inequalities
\begin{equation}\label{eq2-sec2}
|\gamma_{n}|\le \frac{A-B}{2n} \, \mbox{ for } n\ge 1,
\end{equation}
and
\begin{equation}\label{eq1-sec2}
\sum_{n=1}^{\infty}|\gamma_{n}|^{2} \leq \left(\frac{A-B}{2B}\right)^2 {\rm Li\,}_{2}(B^2),
\end{equation}
where ${\rm Li\,}_{2}(x)=\sum_{n=1}^{\infty}\frac{x^n}{n^2}$ denotes the dilogarithm function.
The inequality $(\ref{eq2-sec2})$ is sharp for the function $k_{A,B;n}(z)= z(1+Bz^n)^{(A-B)/nB}$
and the inequality $(\ref{eq1-sec2})$ is sharp for the function $k_{A,B;1}(z)$.
 \end{thm}
%%%%%%%%%%%%%%%%%%%%%%%% Corollary-1 %%%%%%%%%%%%%%%%%%%%%%%%%

%In addition to this case, the following special case will hold.
\begin{cor}\label{cor1}
 If $f\in \es^*(A,-A)$ for $0<A\le 1$, then we have %the sharp inequalities
 $$ |\gamma_n|\le \frac{A}{n} \, \mbox{ for } n\ge 1,\, \mbox{ and } \,
 \sum_{n=1}^{\infty}|\gamma_n|^2 \le {\rm Li_2}(A^2).
 $$
The first and second inequalities are sharp for the functions $k_{A;n}(z)=z/(1-Az^n)^{2/n}$ and $k_{A;1}(z)$, respectively.
\end{cor}

%%%%%%%%%%%%%%%%%%%%%%%%%%%%%%%%%%%%%%%%%%
Theorem \ref{thm1} for the case $B=0$ takes the following form.

\begin{cor}%\label{cor2}
Let $0<A\le 1$ and $f\in\A$ satisfy the inequality $\left|\frac{zf'(z)}{f(z)}-1 \right|<A$, $z\in\D$,
i.e., $f\in\es^*(A,0)$. Then the logarithmic coefficients of $f$ satisfy the inequalities
$$ |\gamma_{n}|\le \frac{A}{2n}\, \mbox{ for } n\ge 1, \mbox{ and }
 \sum_{n=1}^{\infty}|\gamma_{n}|^{2} \leq \frac{A^2}{4}.
$$
Both inequalities are sharp for the functions  $k_{A;n}(z)=ze^{Az^n/n}$ and $k_{A;1}(z)$, respectively.
%{\color{red} Both inequalities are sharp.}
 \end{cor}
 %%%%%%%%%%%%%%%%%%%%%%%% theorem-3 %%%%%%%%%%%%%%%%

Our next result,
which uses the method of proof of \cite[Theorem~2.5]{AvWir-09} and
\cite[Theorem~6.3]{Dur83},
establishes an inequality of the type \eqref{eq:main} for the class $\es^*(A,B)$.

\begin{thm}\label{thm1.2}
Let $f\in\es^*(A,B)$ for $-1\le B<A\le 1$, and let $t\le 2$. Then we have
$$ \sum_{n=1}^{\infty} (n+1)^t|\gamma_n|^2 \le \left(\frac{A-B}{2B}\right)^2\sum_{n=1}^{\infty}\frac{(n+1)^t}{n^2}
|B|^{2n}.
$$
\end{thm}
%%%%%%%%%%%%%%%%%%%%%%%%%%%%

\begin{thm}\label{thm3}
 For $|\alpha|< \pi/2$ and $\beta\in[0,1)$, the logarithmic coefficients of $f\in \es_{\alpha}(\beta)$ satisfy the inequalities
 \begin{equation}\label{eq10-sec2}
 |\gamma_{n}|\le \frac{(1-\beta)}{n}\cos \alpha\, \mbox{ for } n\ge 1,
 \end{equation}
 and
 \begin{equation}\label{eq11-sec2}
   \sum_{n=1}^{\infty}|\gamma_{n}|^{2} \leq \frac{\pi ^2}{6}(1-\beta)^2 \cos^2 \alpha.
 \end{equation}
%  The inequalities are sharp.
Both inequalities are sharp for the function
$f_{\alpha, \beta}(z)=z/(1-z)^{2(1-\beta)\cos \alpha}$. %,~ $\gamma=2(1-\beta)\cos \alpha$.
 \end{thm}
%%%%%%%%%%%%%%%  theorem-3  %%%%%%%%%%%%%%%%%%%%%%%

For the case of strongly starlike functions, we have the following.

\begin{thm}\label{thm4}
 Let $0<\alpha \le1$ and
 $$ A_n(\alpha)=\sum_{k=1}^{n}\binom{n-1}{k-1}\binom{\alpha}{k}2^k.
 $$
 Then the logarithmic coefficients $\gamma_n$ of $f\in \es \es^*_{\alpha}$
 satisfy the inequalities
 \begin{equation}\label{eq1-thm4}
 |\gamma_{n}|\le \frac{\alpha}{n},\quad n\ge 1,
 \end{equation}
 and
 \begin{equation}\label{eq2-thm4}
   \sum_{n=1}^{\infty}|\gamma_{n}|^{2} \leq \frac{1}{4}\sum_{n=1}^{\infty}\frac{|A_n(\alpha)|^2}{n^2}.
 \end{equation}
%where $A_n(\alpha)=\sum_{k=1}^{n}\binom{n-1}{k-1}\binom{\alpha}{k}2^k.$
Both inequalities are sharp for the function $f$ given by $zf'(z)/f(z)= ((1+z^n)/(1-z^n))^{\alpha}$
and $zf'(z)/f(z)= ((1+z)/(1-z))^{\alpha}$, respectively.

 \end{thm}
 %%%%%%%%%%%%%%%%%%%

 Our next two theorems concern the two important classes $ \F(c)$ and  $\mathcal{G}(c)$.

 \begin{thm}\label{thm5}
  Let $f\in \F(c)$ for $c\in (0,3]$. Then the logarithmic coefficients $\gamma_n$ of $f$ for $n=1,2, \ldots, 5$,  satisfy the inequalities
% \begin{equation}\label{eq5*}

$$\left \{\begin{array}{l}
\ds |\gamma_1|\leq \frac{c}{4},\\[2mm]
\ds |\gamma_2|\leq \frac{1}{48}(4c\,+c^2),\\[2mm]
\ds |\gamma_3|\leq \frac{1}{48}(2c\,+c^2),\\[2mm]
% \ds |\gamma_4|\leq \frac{1}{40}\left(c\,+\,\frac{2c^2}{9}\,+\,\frac{c^2}{2}|I_2|\right)\\[2mm]
% \ds |\gamma_5|\leq \frac{1}{60}\left(c\,+\,\frac{c^2}{2}\,+\,\frac{c^3}{24}\,+\,\frac{5c^2}{12}|I_3|\right),
\end{array}
\right .
$$
\[ |\gamma_4|\leq
\left \{\begin{array}{ll}
\ds \frac{1}{40}\left [c\,+\,\frac{c^2}{18}\left(13+\frac{c}{2}-\frac{c^2}{8} \right)\right ] & \mbox{for}~c\in(0,2.61818)\\[4mm]
\ds \frac{1}{40}\left [ c\,+\,\frac{2c^2}{9}\,+\,\frac{c^2}{2}|I_2|\right ] &\mbox{for}~ c\in(2.71569,3]
\end{array}
\right. ,
\]
and
 \[|\gamma_5|\leq
\left \{\begin{array}{ll}
\ds \frac{1}{60}\left [
c\,+\,+\,\frac{c^2}{12}\left(11+c-\frac{c^2}{4} \right)\right ] & \mbox{for}~ c\in(0,1.31148)\\[4mm]
\ds \frac{1}{60}\left [ c\,+\,\frac{c^2}{2}\,+\,\frac{c^3}{24}\,+\,\frac{5c^2}{12}|I_3|\right ] & \mbox{for}~  c\in(1.35541,3]
\end{array}
\right.,
\]
where
\[|I_2|\le \frac{54+c}{27}\left(\frac{4(54+c)}{3(288+8c-c^2)} \right)^{1/2}
\]
and
\[
|I_3|\le \frac{30+c}{15}\left(\frac{2(30+c)}{3(80+4c-c^2)} \right)^{1/2}.
 \]
The first three inequalities are sharp for the functions
$$ f_c(z)
=\left \{\begin{array}{rl}
\ds \frac{(1-z)^{1-c}-1}{c-1} & \mbox{ for $c\neq 1$}\\
-\log (1-z) & \mbox{ for $c =1$.}
\end{array}
\right .
$$
Moreover, the bounds for $|\gamma_4|$ is sharp for $c\le 144/55$ whereas the bounds for $|\gamma_5|$ is sharp for $c\le 80/61$.
\end{thm}
%%%%%%%%%%%%%%% End theorem 5 %%%%%%%%%%%%%%
%%%%%%%%%%%%%%%%%%%%%%%%%%%%%%%%%%%%%%%%%%%%%%%%%%%%%%%%%%%%%%
If we take $c=3$ in Theorem~\ref{thm5},  %by \eqref{eq5*},
then we obtain the logarithmic coefficients bound for the class
$\mathcal{F}(3)$.

\begin{cor}%\label{cor2a}
Let $\mathcal{F}(3)$, i.e. ${\rm Re}  \left ( 1+(zf''(z)/f'(z))\right )>-1/2$ in $\D$. Then we have
\[ |\gamma_1|\leq \frac{3}{4},\,\,|\gamma_2|\leq \frac{7}{16},~~ |\gamma_3|\leq \frac{5}{16},\,\]
\[ |\gamma_4|\leq \frac{1}{40}\left(5\,+\,\frac{19}{2}\sqrt{\frac{76}{303}}\right)\approx 0.243945
\]
and
\[ |\gamma_5|\leq \frac{1}{60}\left(\frac{69}{8}\,
+\,\frac{33}{4}\sqrt{\frac{22}{83}}\right)\approx 0.2145050.
\]
The first three inequalities are sharp for the function
$${f_3(z)= \frac{z-(z^2/2)}{(1-z)^2}}
$$
which is extremal for the class $\mathcal{F}(3)$. For this function, we have
\begin{equation}\label{eq-cor2}
\gamma_k(f_3)=\left(1-\frac{1}{2^{k+1}} \right)\frac{1}{k} ~\mbox{  for $k\in \mathbb{N}$}.
\end{equation}
\end{cor}
%%%%%%%%%%%%%%%%%%%%%%%%%%%%%%%%%%%%%%%%%%%%%%%%%%%

The logarithmic coefficients  of $f_3$ given by \eqref{eq-cor2} gives that
\begin{align*}
 4\sum_{n=1}^{\infty } |\gamma_n(f_3)|^2 &  =
  \sum_{n=1}^{\infty} \left(2-\frac{1}{2^n} \right)^2 \frac{1}{n^2}\\
 & = 4\sum_{n=1}^{\infty} \frac{1}{n^2}+\sum_{n=1}^{\infty}\frac{(1/4)^n}{n^2}-
 4\sum_{n=1}^{\infty}\frac{(1/2)^n}{n^2}\\
 & = \frac{4\pi^2}{6}+ {\rm Li\,}_{2}\left(\frac{1}{4}\right)
    -4~{\rm Li\,}_{2}\left(\frac{1}{2}\right),
 \end{align*}
where ${\rm Li\,}_{2}(x)=\sum_{n=1}^{\infty}\frac{x^n}{n^2}$ denotes the dilogarithm function.
These observations lead us to the following conjecture.
\bcon
The logarithmic coefficients $\gamma_n$ of $f\in\F (3)$ satisfy the inequalities
$$ |\gamma_n|\le \frac{1}{n}\left(1-\frac{1}{2^{n+1}} \right) \, \mbox{ for } n\ge 1
$$
and
$$ \sum_{n=1}^{\infty}|\gamma_{n}|^{2} \leq \frac{\pi^2}{6}+\frac{1}{4} ~{\rm Li\,}_{2}\left(\frac{1}{4}\right)
   -{\rm Li\,}_{2}\left(\frac{1}{2}\right).
$$
Equalities in these inequalities are attained for the function $f_3(z)$ given as above.
% Both inequalities are sharp for the function $f_3(z)$.
\econ

In \cite{OPW18}, the authors considered the logarithmic coefficients of the functions $f$ in the class
$\mathcal{G}(c)$
% that satisfy the inequality $$ \R\left(1\,+\,\frac{zf''(z)}{f'(z)}\right)\,<1\,\,+\,\frac{c}{2},\quad z\in \D$$
for some $c\in (0,1]$ and they got the estimate
$$|\gamma_n|\leq \frac{c}{2(c+1)n},\quad n\in \N.$$
Among others, they conjectured that for $c=1$ the inequalities
$$|\gamma_n|\leq \frac{1}{2(n+1)n},\quad n\in \N,
$$
are valid, where equality is attained for $f'(z)\,=\,(1-z^n)^{1/n}.$ In \cite{OPW18}, this conjecture was proved only in the case $n=1$.
In the sequel, we shall consider the cases $n=1,2,3$ in the families $\mathcal{G}(c)$ using a
similar method as in the proof of Theorem~\ref{thm5}.

%%%%%%%%%%%%%%%%%% Theorem-6 for G(c) %%%%%%%%%%%%%%%
\begin{thm}\label{thm6}
Let $f\in \G(c)$ for $c\in (0,1]$. Then the logarithmic coefficients $\gamma_n$ of $f$  satisfy the inequalities
% \begin{equation}\label{eq5*}
\[ |\gamma_1|\leq \frac{c}{4},\, |\gamma_2|\leq \frac{c}{12},\, \mbox{ and } \,
|\gamma_3|\leq \frac{c}{24}.
\]
The inequalities are attained for $f'(z)\,=\,(1-z^n)^{c/n},\,n=1,2,3.$
\end{thm}
%%%%%%%%%%%%%%%%%%%%%%
% \bigskip

These results led us to a generalization of the conjecture mentioned above:

\bcon
The logarithmic coefficients $\gamma_n$ of the functions in $\mathcal{G}(c),\, c\in (0,1],$ satisfy the inequalities
$$|\gamma_n|\,\leq \,\frac{c}{2n(n+1)},\quad n\in \N.
$$
Equality is attained for $f'(z)\,=\,(1-z^n)^{c/n}.$
\econ

%%%%%%%%%%%%%%%%%%%%%%%
%%%%%%%%%%%%%%%%%%%%%%%%%%%%%%%%%%
\section{Lemmas}\label{lemmas}

\begin{lem}\cite[Corollary~3.1d.1, p.~76]{MM2000}\label{lem1}
 Let $h$ be starlike in $\D$, with $h(0)=0$ and $a\neq 0$. If an analytic function
 $p(z)=a+a_nz^n+a_{n+1}z^{n+1}+\cdots$ satisfies the subordination relation
 $$ \frac{zp'(z)}{p(z)} \prec h(z),
 $$
 then
 $$ p(z) \prec q(z)=a\exp\left[ n^{-1}\int_0^z\frac{h(t)}{t}~ dt\right].
 $$
\end{lem}
%%%%%%%%%%%%%%%%%%%%%%%%%%%%%%%%%
\begin{lem}\cite[Theorem~6.3, p.~192]{Dur83}~
$($see also \cite[Rogosinski's Theorem~II (i)]{Rog43}$)$\label{lem2}
Let $f(z)=\sum_{n=1}^{\infty}a_nz^n$ and $g(z)=\sum_{n=1}^{\infty}b_nz^n$ be analytic in $\D$,
and suppose that $f \prec g$, where $g$ is univalent in $\D$. Then
$$ \sum_{k=1}^{n}|a_k|^2 \le \sum_{k=1}^{n}|b_k|^2,\quad n=1,2,\dots.
$$
\end{lem}
%%%%%%%%%%%%%%%%%%%%%%%%%
%%%%%%%%%%%%%%%%%%%%%%%%%%%%
\begin{lem}\cite[Theorem~6.4 (i), p.~195]{Dur83}~$($see also \cite[Rogosinski's Theorem~X]{Rog43}$)$\label{lem4}
Let $f(z)=\sum_{n=1}^{\infty}a_nz^n$ and $g(z)=\sum_{n=1}^{\infty}b_nz^n$ be analytic in $\D$. Suppose that
$f \prec g$, where $g$ is univalent and convex in $\D$. Then
$$ |a_n|\le |g'(0)|=|b_1|,\quad n=1,2,\dots.
$$
\end{lem}
%%%%%%%%%%%%%%%%%%%%%%%%%%%

Our next lemma due to Prokhorov and  Szynal \cite{PS} is crucial in the investigation of fourth and fifth
logarithmic coefficients bound for $\F(c)$ and $\G(c)$.

\begin{lem}\cite[Lemma~2]{PS}\label{lem-7}
Let $\phi(z)= \sum_{k=1}^{\infty}c_kz^k \in \mathcal B$ be a Schwarz function and
$$\Psi(\phi)=|c_3+\mu c_1c_2+\upsilon c_1^3|.
$$
Then for any real number $\mu$ and $\upsilon$, we have the following sharp estimate:
\[
\Psi(\phi)\le \Phi(\mu,\upsilon)=
\left\{
\begin{array}{ccc}
1 && \mbox{if }(\mu,\upsilon)\in D_2,\\[1mm]
|\upsilon| &&\mbox{if }(\mu,\upsilon)\in D_6,\\[1mm]
\displaystyle{ \frac{2}{3}\,(|\mu|+1)\left(\frac{|\mu|+1}{3\left(|\mu|+1+\upsilon\right)}\right)^{1/2} }
       && \mbox{if }(\mu,\upsilon)\in D_9,
\end{array}
 \right.
 \]
where
\begin{align*}
 D_2 &=\left\{(\mu,\upsilon)\in\mathbb{R}^2:\, \frac{1}{2}\le |\mu|\le 2,~~ \frac{4}{27}(|\mu|+1)^3-(|\mu|+1) \le \upsilon \,\le 1 \right\},\\[1mm]
 D_6 &= \left\{(\mu,\upsilon)\in\mathbb{R}^2:\, 2\le |\mu|\le 4,~~  \upsilon \,\ge \frac{1}{12}(\mu^2+8)\right\}, ~\mbox{ and }\\[1mm]
 D_9 &=\left\{(\mu,\upsilon)\in\mathbb{R}^2:\, |\mu|\ge 2,~~ -\frac{2}{3}\,(|\mu|+1)\,\le \upsilon \,\le
\frac{2|\mu|\,(|\mu|+1)}{\mu^2+2|\mu|+4}  \right\}.
 \end{align*}
\end{lem}
%%%%%%%%%%%%%%%%%%%%%%%%%%%%%%%%%%%%%%%%%%%
\section{Proofs of the main results} \label{proof theorems}

\subsection{Proof of Theorem~\ref{thm1}}
Suppose $f \in \es^*(A,B)$ for  $-1\leq B<A\le 1$. Then by the definition of $\es^*(A,B)$, we get
\begin{equation}\label{eq5-sec2}
 z\left (\log\left (\frac{f(z)}{z}\right )\right )'=
 \frac{zf'(z)}{f(z)}-1\prec %\frac{1+Az}{1+B z}-1 =
 \frac{(A-B)z}{1+B z}, \quad z\in \D,
\end{equation}
which, in terms of the logarithmic coefficients $\gamma_{n}$ of $f$ defined by (\ref{eq1-sec1}), is equivalent to
\begin{align}\label{eq6-sec2}
 \sum_{n=1}^{\infty}(2n \gamma_{n})z^{n}\prec &
 \left \{
\begin{array}{ll}
\ds  \nonumber \frac{(A-B)}{B} \sum_{n=1}^{\infty}(-1)^{n-1}B^nz^n
 & \mbox{ if } -1\le B<A\le 1, B\neq 0\\[4mm]
 Az & \mbox{ if } B= 0
 \end{array}\right.\\
 & =:G(z).
\end{align}
Since $G$ is convex in $\D$ with $G'(0)=A-B$, it follows from
%Rogosinski's result
% (see also \cite[Theorem 6.4(i), p.~195]{Dur83} or see Lemma~\ref{lem4})
Lemma~\ref{lem4} that
$$ 2n|\gamma_n|\le |G'(0)|=|A-B|\, \mbox{ for } n\ge 1,
$$
which implies the desired inequality (\ref{eq2-sec2}).
The equality holds for the function $k_{A,B;n}(z)=  z(1+Bz^n)^{(A-B)/nB}$.
We have
$$ \log \left (\frac{k_{A,B;n}(z)}{z}\right ) =\frac{A-B}{nB}\sum_{k=1}^{\infty}\frac{(-1)^{k+1}B^k}{k}(z^{n})^k. %=:2\sum_{n=1}^{\infty}\gamma_n(k_n)z^n,
$$
%where $\gamma_n(k_n)=(A-B)/2n.$

Let $g(z):=z/f(z)$ which is a non-vanishing analytic function in $\D$ with $g(0)=1$, and it has the series representation
$$ g(z)=1+\sum_{n=1}^{\infty}b_nz^n.
$$
It is clear from (\ref{eq5-sec2}) that $g$ satisfies the relation
$$ \frac{zg'(z)}{g(z)}=1-\frac{zf'(z)}{f(z)}\prec \frac{-(A-B)z}{1+B z}=:\phi(z), \quad z\in \D.
$$
Note that $\phi$ is convex in $\D$ and $\phi(0)=0$. By using the subordination result of Lemma~\ref{lem1}, we get
\begin{equation}\label{eq9-sec2}
g(z):=\frac{z}{f(z)} \prec q_{A,B}(z)=\exp \left(\int_0^{z}\frac{\phi(t)}{t} dt \right).
%q_{A,B}(z)=
%\left \{
%\begin{array}{ll}
%(1+Bz)^{1-(A/B)} & \mbox{ for } B\neq 0\\
%e^{-Az} & \mbox{ for } B=0.
%\end{array} \right.
\end{equation}
It is a simple exercise to compute that
$$ \displaystyle q_{A,B}(z)=
 \left \{
 \begin{array}{ll}
 e^{-Az} & \mbox{ for } B=0\\
 (1+Bz)^{1-(A/B)} & \mbox{ for } B\neq 0.
 \end{array} \right.
$$
We can rewrite the relation (\ref{eq9-sec2}) as
$$ \frac{f(z)}{z}\prec \frac{1}{q_{A,B}(z)},
$$
which, in terms of the logarithmic coefficients $\gamma_{n}$ of $f$ defined by (\ref{eq1-sec1}), is equivalent to
(compare with (\ref{eq6-sec2}))
$$ 2\sum_{n=1}^{\infty} \gamma_{n}z^{n} \prec
\left \{
\begin{array}{ll}
\ds
\left(\frac{A-B}{B}\right)\sum_{n=1}^{\infty} (-1)^{n-1}\frac{B^n}{n}z^{n}
= \left(\frac{A-B}{B}\right) \log (1+Bz)\\[4mm]
\hspace{5.5cm}\mbox{ for } B\neq 0\\[4mm]
Az
\hspace{5cm} \mbox{ for } B= 0.
\end{array}\right.
$$
% for $B\neq 0$.
% By Rogosinski's result (see \cite[Theorem~II (i)]{Rog43} or see Lemma~\ref{lem2}),
Using Lemma~\ref{lem2}, we obtain that
\begin{align*}
 4\sum_{n=1}^{k} |\gamma_{n}|^2 & \le \left(\frac{A-B}{B}\right)^2\sum_{n=1}^{k} \frac{B^{2n}}{n^2}
 \quad \mbox{ for }B\neq 0\\
  & \le (A-B)^2 \frac{{\rm Li\,}_{2}(B^2)}{B^2}.
\end{align*}
Letting $k \rightarrow \infty$, we get
$$ \sum_{n=1}^{\infty} |\gamma_{n}|^2 \le \frac{(A-B)^2}{4} \frac{{\rm Li\,}_{2}(B^2)}{B^2},
$$
where ${\rm Li\,}_{2}(x)=\sum_{n=1}^{\infty}\frac{x^n}{n^2}$. For $x=0$, we let $\frac{{\rm Li\,}_{2}(x)}{x}$ as the limit value 1.
This proves the desired assertion (\ref{eq1-sec2}).
The equality holds for the function $k_{A,B}$  defined by
$$ \displaystyle k_{A,B}(z):=
\left \{
\begin{array}{ll}
ze^{Az} & \mbox{ for } B=0\\
z(1+Bz)^{(A/B)-1} & \mbox{ for } B\neq 0.
\end{array} \right.
$$
%is in $\es ^*(A,B)$ and acts the role of extremal function for this class.

Indeed, for the function $k_{A,B}$, we have
\begin{align*}
\displaystyle \log \left (\frac{k_{A,B}(z)}{z}\right ) & =
\left \{
\begin{array}{ll}
 \ds \left (\frac{A-B}{B}\right )\log (1+Bz) & \mbox{ for } B\neq 0\\[2mm]
{Az} & \mbox{ for } B=0
 \end{array} \right. \\[3mm]
 & = 2\sum_{n=1}^{\infty} \gamma_{n}(k_{A,B})z^{n},
\end{align*}
where
$$
\displaystyle \gamma_{n}(k_{A,B}) =
\left \{
\begin{array}{ll}
\ds
(-1)^{n-1} \left (\frac{A-B}{2B}\right )\frac{B^n}{n}& \mbox{ for } B\neq 0\\[2mm]
\ds \frac{A}{2} & \mbox{ for } B=0.
 \end{array} \right. \\
$$
This completes the proof of Theorem~\ref{thm1}.
\hfill$\Box$

\vspace{.5cm}
We remark that when $A=1-2\beta$ and $B=-1$ in Theorem~\ref{thm1}, we obtain~\cite[Remark~1]{OPW18}.
%%%%%%%%%%%%%%%%%%%%%%%%%%%%%%%%

%%%%%%%%%%%% remark-1 %%%%%%%%%%
\begin{rem}\label{rem1}
 From the relation $(\ref{eq6-sec2})$ and using Rogosinski's theorem $($see \cite[Theorem~6.3]{Dur83}$)$, we obtain
 $$ 4\sum_{n=1}^{\infty}n^2|\gamma_n|^2 \le \left|\frac{A-B}{B}\right|^2\sum_{n=1}^{\infty}|B|^{2n}\, \mbox{ for } B\neq -1
 $$
 and so
 \begin{equation}\label{eq1-rem1}
  \sum_{n=1}^{\infty}n^2|\gamma_n|^2 \le \frac{(A-B)^2}{4(1-B^2)}.
 \end{equation}
\end{rem}
%%%%%%%%%%%%%%%%%%%%%
\subsection{Proof of Theorem~\ref{thm1.2}}
We recall from the formula $(\ref{eq6-sec2})$ and the result of Rogosinski
$($see also~\cite[Theorem 2.2]{Pom75} and~\cite[Theorem 6.3]{Dur83}$)$, that for $k\in \N$ the inequalities
\begin{equation}\label{eq3-rem1}
4\sum_{n=1}^{k} n^2|\gamma_n|^2 \le \left(\frac{A-B}{B}\right)^2\sum_{n=1}^{k}|B|^{2n}
\end{equation}
are valid. This implies the inequality $(\ref{eq1-rem1}$) as well, if  $B\neq -1$.
We consider $(\ref{eq3-rem1})$ for $k=1,2,\dots, N$, and
multiply the $k$-th inequality by the factor
$$\frac{(k+1)^t}{k^2}-\frac{(k+2)^t}{(k+1)^2}>0, \mbox{ if } k=1,2,\dots, N-1
$$
and by $\frac{(N+1)^t}{N^2}$ for $k=N$.
Then the summation of the multiplied inequalities yields
% Adding up these modified inequalities results in the inequality
% \begin{equation}\label{eq1-rem2}
$$ \sum_{k=1}^{N} (k+1)^t|\gamma_k|^2  \le \left(\frac{A-B}{2B}\right)^2\sum_{k=1}^{N}
\frac{(k+1)^t}{k^2}B^{2k}
% &  \le \left(\frac{A-B}{2B}\right)^2\sum_{k=1}^{\infty}\frac{1}{k^2}(B)^{2k}\\
% & =  \left(\frac{A-B}{2B}\right)^2 {\rm Li_2}(B^2),
% \end{equation}
$$
for $N=1,2,\dots$. Allowing $N\rightarrow \infty$, we see that the proof of the theorem is complete.
\hfill$\Box$
%%%%%%%%%%%%%%%%%%%%%% remark-2 %%%%%%%%%%%%%%%
\begin{rem}\label{rem2}
Here is an alternate approach to prove the inequality ~$(\ref{eq1-sec2})$
%According to this method of proof, compare for example \cite[Theorem~2.5]{AvWir-09} and \cite[Theorem~6.3]{Dur83}.
(compare for example \cite[Theorem~2.5]{AvWir-09}).
If we take $t=0$ in Theorem~\ref{thm1.2}, then we obtain
$$
\sum_{n=1}^{\infty} |\gamma_n|^2
\le \left(\frac{A-B}{2B}\right)^2\sum_{n=1}^{\infty}\frac{1}{n^2}
B^{2n} = \left(\frac{A-B}{2B}\right)^2 {\rm Li_2}(B^2).
$$
% for $N=1,2,\dots$. Allowing $N\rightarrow \infty$, we get $(\ref{eq1-sec2})$ again.
\end{rem}
%%%%%%%%%%%%%%%%%%%%%%%%
%%%%%%%%%%  End proof of Theorem 1   %%%%%%%%%%%%%%%%%%%%%%%%%%%%%%%%%%%%%%%
%%%%%%% Proof of Theorem-3 %%%%%%%%%%%
\subsection{Proof of Theorem~\ref{thm3}}
% \begin{proof}
 Suppose $f\in\es_{\alpha}(\beta)$. Then by the definition of $\es_{\alpha}(\beta)$, we obtain that
\beq\label{eq12-sec2}
 \frac{zf'(z)}{f(z)}-1  \prec  & \ds \frac{1+[e^{2i\alpha}-2\beta e^{i\alpha}\cos \alpha]z}{1-z}-1 \nonumber \\
 & \ds = 2(1-\beta)e^{i\alpha}\cos \alpha\left(\frac{z}{1-z}\right)=:G_1(z),
\eeq
which, in terms of the logarithmic coefficients $\gamma_{n}$ of $f$ defined by (\ref{eq1-sec1}), is equivalent to
$$ 2\sum_{n=1}^{\infty}n \gamma_{n}z^{n}\prec G_1(z).
$$
Since $G_1$ is convex in $\D$ and $G_1'(0)=2(1-\beta)e^{i\alpha}\cos \alpha$,
Rogosinski's result (see Lemma~\ref{lem4}) gives the inequality (\ref{eq10-sec2}).

From (\ref{eq12-sec2}), $g$ defined by $g(z)=z/f(z)$ satisfies the relation
$$ \frac{zg'(z)}{g(z)}\prec -G_1(z), \quad z\in \D.
$$
Note that $G_1$ is convex in $\D$ with $G_1(0)=0$. By Lemma~\ref{lem1}, we get
$$ g(z)=\frac{z}{f(z)} \prec q_{\alpha,\beta}(z):=\exp\left(-\int_0^z \frac{G_1(t)}{t}~dt\right),
$$
or equivalently,
$$ \frac{f(z)}{z} \prec \frac{1}{q_{\alpha,\beta}(z)}= (1-z)^{-\gamma}=:\frac{f_{\alpha,\beta}(z)}{z},
$$
where $\gamma= 2(1-\beta)\cos \alpha$. It is easy to see that $f_{\alpha,\beta} \in \es_{\alpha}(\beta)$.
From (\ref{eq1-sec1}), we have
$$ 2\sum_{n=1}^{\infty}\gamma_{n}z^{n} \prec \gamma \sum_{n=1}^{\infty} \frac{z^{n}}{n}.
$$
Then by  %Rogosinski's theorem (see Lemma~\ref{lem2}),
Lemma~\ref{lem2}, we obtain
$$ \sum_{n=1}^{k}|\gamma_{n}|^2 \le \frac{|\gamma|^2}{4}\sum_{n=1}^{k} \frac{1}{n^2} \le
\frac{|\gamma|^2}{4}\sum_{n=1}^{\infty} \frac{1}{n^2}.
$$
If we allow $k\rightarrow \infty$, we get the inequality (\ref{eq11-sec2}).
The equality holds in the inequalities (\ref{eq10-sec2}) and (\ref{eq11-sec2}) for the function
$f_{\alpha, \beta}(z)=z/(1-z)^{\gamma}$,~ $\gamma=2(1-\beta)\cos \alpha$.
% \end{proof}
% The case $\alpha=0$ in Theorem~\\hfill$\Box$ref{thm3} is recently obtained in~\cite[Remark~1]{OPW18}.
\hfill$\Box$

%%%%%%%%%%%%%%%%%%%% for the class of Strongly starlike function %%%%%%%%%%%%%%%%%%%%
\subsection{Proof of Theorem~\ref{thm4}}
% \begin{proof}
 Suppose $f\in \es \es^*_{\alpha}$, and $g(z)=z/f(z)$. Then $g$ is a non-vanishing analytic function in $\D$
 and by the definition of $\es \es^*_{\alpha}$, we get
\begin{equation}\label{eq3-thm4}
 \frac{zg'(z)}{g(z)}=1-\frac{zf'(z)}{f(z)} \prec 1-\left(\frac{1+z}{1-z} \right)^{\alpha}=:\phi_1(z),\quad z\in\D.
\end{equation}
Using Lemma~\ref{lem1}, we obtain
$$ g(z)=\frac{z}{f(z)} \prec q_{\alpha}(z):=\exp\left(\int_0^z \frac{\phi_1(t)}{t}\right)~dt,
$$
which is equivalent to
$$ \log \left (\frac{f(z)}{z}\right ) \prec \log \left(\frac{1}{q_{\alpha}(z)}\right)=-\int_0^z \frac{\phi_1(t)}{t}~dt.
$$
We consider the function
\begin{align}\label{eq5-thm4}
\left(\frac{1+z}{1-z} \right)^{\alpha} = 1+ \sum_{n=1}^{\infty} A_n(\alpha)z^n.
\end{align}
Using this, we have
$$ \int_0^z \frac{\phi_1(t)}{t}~ dt = -\sum_{n=1}^{\infty} A_n(\alpha)\frac{z^{n}}{n},
$$
where
$$A_n(\alpha)=\sum_{k=1}^{n}\binom{n-1}{k-1}\binom{\alpha}{k}2^k\quad \mbox{ for } n\ge 1.
$$
 As in the previous case, by Rogosinski's theorem (see Lemma~\ref{lem2}), we obtain
$$ \sum_{n=1}^{\infty}|\gamma_{n}|^2 \le  \frac{1}{4}\sum_{n=1}^{\infty} \frac{|A_n(\alpha)|^2}{n^2}.
$$
The proof of (\ref{eq2-thm4}) is completed.

Now, we prove the inequality (\ref{eq1-thm4}). From (\ref{eq3-thm4}) and (\ref{eq5-thm4}), we also get
$$ \frac{zf'(z)}{f(z)} -1 \prec -\phi_1(z)=\sum_{n=1}^{\infty} A_n(\alpha)z^n,
$$
which, in terms of the logarithmic coefficients $\gamma_{n}$ of $f$ defined by (\ref{eq1-sec1}), is equivalent to
$$2\sum_{n=1}^{\infty} n \gamma_{n}z^{n} \prec \sum_{n=1}^{\infty} A_n(\alpha)z^n.
$$
Using Lemma~\ref{lem4}, we find that
$$  2n |\gamma_{n}| \le |-\phi_1'(0)|=|A_1(\alpha)|=2\alpha,
$$
and the desired inequality (\ref{eq1-thm4}) follows.
Equality occurs in the inequalities (\ref{eq1-thm4}) and (\ref{eq2-thm4})
if $f\in\A$ is given by
$$\frac{zf'(z)}{f(z)} = \left (\frac{1+z^n}{1-z^n}\right )^{\alpha}
~\mbox{ and } ~\frac{zf'(z)}{f(z)}=  \left (\frac{1+z}{1-z}\right )^{\alpha},
$$
respectively.
% \end{proof}
\hfill$\Box$

% We remark that when $\alpha=1$ in Theorem~\ref{thm4}, we obtain \cite[Remark~1]{OPW18}.
%%%%%%%%%%%%% End Theorem 4 %%%%%%%%%%%%%%%%%%
\subsection{Proof of Theorem~\ref{thm5}}
Let $f\in \F(c)$ for $c\in (0,3]$. Consider the identity
\begin{equation}\label{eq1}
\left(\frac{zf'(z)}{f(z)}\right) \left(1\,+\,\frac{zf''(z)}{f'(z)}\right)\,
=\,\left(\frac{zf'(z)}{f(z)}\right)^2\,
+\,z\left(\frac{zf'(z)}{f(z)}\right)'.
\end{equation}
%  and consider
% $$-\frac{1}{2}\log\left (\frac{z}{f(z)}\right )=\sum_{n=1}^\infty \gamma_nz^n.
% $$
Now, we may set
\begin{equation}\label{eq1.1}
1\,+\,\frac{zf''(z)}{f'(z)} =\,\sum_{n=0}^{\infty}\beta_nz^n,
\end{equation}
where $\beta_0\,=\,1$. %By the previous relation, we have , (\ref{eq1}) and (\ref{eq1.1})
By the relation (\ref{eq1-sec1}), we have
%\[\frac{zf'(z)}{f(z)}-1\,=\,2\sum_{n=1}^{\infty}n\gamma_nz^n=\,\sum_{n=1}^{\infty}\delta_nz^n,
%\quad \delta_n=2n\gamma_n,\]
%or
\begin{equation}\label{eq2}
\frac{zf'(z)}{f(z)}=\,1+2\sum_{n=1}^{\infty}n\gamma_nz^n=\,\delta_0+\sum_{n=1}^{\infty}\delta_nz^n,
\end{equation}
where $\delta_0\,=\,1$ and $\delta_n=2n\gamma_n$.
% From (\ref{eq1}), using (\ref{eq1.1}) and (\ref{eq2}),
%%%%%%%%%%%%%%%%%%%%%%%%%%%%%
% From (\ref{eq1}) and (\ref{eq2}),
Using (\ref{eq1.1}) and (\ref{eq2}), we can write (\ref{eq1}) in the series form  as
\[ \left(\sum_{n=0}^{\infty}\beta_n z^n \right) \left(\sum_{n=0}^{\infty}\delta_n z^n \right)
= \left(\sum_{n=0}^{\infty}\delta_n z^n \right)^2 + \sum_{n=0}^{\infty}n\delta_n z^n.
 \]
Using the Cauchy product of power series, we obtain
\begin{align*}
\sum_{n=0}^{\infty}\left(\sum_{k=0}^{n}\delta_n \beta_{n-k} \right)z^n & =
\delta_0\beta_0+(\delta_0\beta_1+\delta_1\beta_0)z + (\delta_0\beta_2+\delta_1\beta_1+\delta_2\beta_0)z^2+\cdots\\
& = (\delta_0+\delta_1 z+\delta_2 z^2+\delta_3 z^3+\cdots)^2+ \delta_1z+2\delta_2z^2\\
&\hspace{1cm} +3\delta_3z^3+\cdots.
\end{align*}
As $\beta_0=1=\delta_0$, it is equivalent to
\begin{align}\label{eq2.1}
& 1+(\delta_1+\beta_1)z+(\delta_2+\delta_1\beta_1+\beta_2)z^2+(\beta_3+\delta_1\beta_2+\delta_2\beta_1+\delta_3)z^3
+\cdots \nonumber\\
& \hspace{.2cm}= 1+(2\delta_1+\delta_1)z+(2\delta_2+\delta_1^2+2\delta_2)z^2
+(2\delta_3+2\delta_1\delta_2+3\delta_3)z^3+\cdots.
\end{align}
%%%%%%%%%%%%%%%%%%%%%%%%%%%%%%%
First, we compare the coefficients of $z^n$ for $n=1,2,3$ and get (by using $\gamma_n=\delta_n/(2n)$) that
\begin{equation}\label{eq3}
\gamma_1\,=\,\frac{\delta_1}{2}\,=\,\frac{\beta_1}{4},
\end{equation}
\begin{equation}\label{eq4}
\gamma_2\,=\,\frac{\delta_2}{4}\,=\,\frac{1}{12}\left(\beta_2\,+\,\frac{\beta_1^2}{4}\right),
\end{equation}
and
\begin{equation}\label{eq5}
\gamma_3\,=\,\frac{\delta_3}{6}\,=\,\frac{1}{24}\left(\beta_3\,+\,\frac{\beta_1\beta_2}{2}\right).
\end{equation}
% Each $f\in \mathcal{F}(c)$ can be written equivalently, in terms of subordination and \eqref{eq1.1}, as
By the definition of $f\in \F (c)$  and \eqref{eq1.1}, we get
\begin{equation}\label{eq5.11}
\sum_{k=0}^{\infty}\beta_kz^k \, =\,
1\,+\,\frac{zf''(z)}{f'(z)}\,\prec\,1\,+\,\frac{cz}{1-z},\quad z\in \mathbb{D}.
\end{equation}
But then,  by Rogosinski's Theorem (see, for instance, \cite[p.~195, Theorem 6.4]{Dur83}), we obtain that
% \begin{equation}\label{eq5.1}
$$ |\beta_k|\leq c,\quad k\in \mathbb{N}$$
% \end{equation}
which also  follows from the estimate on the coefficients of functions with real part bigger than
$\alpha=1-c/2 < 1:$
$$|\beta_k|\leq 2(1-\alpha)=c, \quad k\in \mathbb{N}.
$$
The estimate here is sharp. Using this inequality, the above identities \eqref{eq3}--\eqref{eq5} imply
% \begin{equation}\label{eq5*}
$$ |\gamma_1|\leq \frac{c}{4},~|\gamma_2|\leq \frac{1}{48}(4c\,+c^2), ~\mbox{ and }~
|\gamma_3|\leq \frac{1}{48}(2c\,+c^2).
$$
% \end{equation}
These inequalities are sharp and the equality follows from \eqref{eq3}--\eqref{eq5} by substituting $\beta_k=c$ for $k=1,2,3$.
The extremal function can be derived by integration of the equation
\[1\,+\,\frac{zf''(z)}{f'(z)}\,=\,1\,+\,\frac{cz}{1-z},\]
which results in
\[
f(z)\,=\,\frac{(1-z)^{1-c}-1}{c-1} =:f_c(z) \mbox{ for $c\neq 1$}\]
and $f_1(z)\,=\,-\log(1-z)$ for  $c=1.$

%%%%%%%%%%%%%%%%%%%%%%%%%%%%%%%%%
To get estimates for the fourth and fifth logarithmic coefficients, by (\ref{eq5.11}), we use the fact that
\begin{equation}\label{eq-5.12}
\sum_{k=1}^{\infty}\beta_kz^k\,=\,\frac{c\phi(z)}{1-\phi(z)},
\end{equation}
where $\phi \in\mathcal{B}$ is a Schwarz function. Note that equality is attained in the estimate
$|\beta_k|\leq c$, if $\phi(z)\,=\,e^{i\theta}z.$
If we write $\phi(z)\,=\,\sum_{k=1}^{\infty}c_kz^k,$ then the identity (\ref{eq-5.12}) implies
\[
\beta_1\,=\,cc_1,\quad \beta_2\,=\,c\left(c_2\,+\,c_1^2\right),~\mbox{ and }~
\beta_3\,=\,c\left(c_3\,+\,2c_1c_2\,+\,c_1^3\right).\]
In view of these, expressions of the form
\[
I_1=\beta_3\,+\,a\beta_1\beta_2\,+\,b\beta_1^3\]
can be written in the form
\begin{equation}\label{eq-s1}
I_1=c\left[c_3\,+\,(2+ac)c_1c_2\,+\,c_1^3(1+ac+bc^2)\right].
\end{equation}
%%%%%%%%%%%%% New Lemma %%%%%%%%%%%%%%%%%%%%%%%%%%%%
We now consider the functional
\[ \Psi(\phi)=|c_3+\mu c_1c_2+\upsilon c_1^3|,\, \mbox{ $\mu$ and $\upsilon$ are real},
\]
within the class $\mathcal B$.
In \cite{PS}, Prokhorov and Szynal found the precise bound for the functional $\Psi(\phi)$
which we recalled in Lemma~\ref{lem-7} and use this to find the fourth and the fifth logarithmic
coefficients bound for the class $\mathcal{F}(c)$.

%%%%%%%%%%%%%%%%%%%%%%%%%%%%%%%%%%%%%%%%%%%
Now, we calculate the fourth logarithmic coefficient bound for the class $\mathcal{F}(c)$.
From (\ref{eq2.1}), we compare the coefficients of $z^4$ and obtain
\[ \beta_4+\beta_3\delta_1+\beta_2\delta_2+\beta_1\delta_3+\delta_4
= 6\delta_4+2\delta_1\delta_3+\delta_2^2.
 \]
Using (\ref{eq3}), (\ref{eq4}) and then simplifying, we get
\begin{align*}
 5\delta_4 & = \beta_4+\beta_3\delta_1+\beta_2\delta_2- \delta_2^2
             = \beta_4+\beta_3\delta_1+\delta_2(\beta_2- \delta_2)\\
           & = \beta_4+ \frac{\beta_3\beta_1}{2}+ \frac{1}{3}\left(\beta_2\,+\,\frac{\beta_1^2}{4}\right)
            \left[\beta_2-\frac{1}{3}\left(\beta_2\,+\,\frac{\beta_1^2}{4}\right)\right].
\end{align*}
Since $\delta_4=8\gamma_4$, we have
\begin{align}\label{eq5.2}
 \nonumber
 40 |\gamma_4| & = \left|\beta_4+ \frac{\beta_3\beta_1}{2}+ \frac{1}{9}\left(2\beta_2^2\,
                    +\,\frac{\beta_1^2 \beta_2}{4}-\,\frac{\beta_1^4}{16}\right)\right|\\ \nonumber
               & \le |\beta_4|\, +\, \frac{2}{9}|\beta_2|^2 \,+\, \frac{|\beta_1|}{2}\left|\beta_3 + \,
                    \frac{\beta_1 \beta_2}{18} -\,\frac{\beta_1^3}{72}\right|\\ \nonumber
               & \le c\,+\,\frac{2c^2}{9}\,+\,\frac{c^2}{2}\left|c_3\,+\left(2+\frac{c}{18}\right)c_1c_2\,
                    +\,\left(1+\frac{c}{18}-\frac{c^2}{72}\right)c_1^3\right| \\ \nonumber
                    &     \qquad (\mbox{using (\ref{eq-s1})})\\
               &   =\,c\,+\,\frac{2c^2}{9}\,+\,\frac{c^2}{2}|I_2|,
\end{align}
where
\begin{equation}\label{eq1ex}
I_2
%= c_3\,+\left(2+\frac{c}{18}\right)c_1c_2\,+\,\left(1+\frac{c}{18}-\frac{c^2}{72}\right)c_1^3\\
:= c_3\,+\mu c_1c_2\,+\upsilon c_1^3, \quad \mu =2+\frac{c}{18} ~\mbox{ and }~\upsilon =1+\frac{c}{18}-\frac{c^2}{72}.
\end{equation}
Since our aim is to get an upper bound that corresponds to $|c_1|=1$, we ask for values of $c$ such that
\begin{equation}\label{eqa}
% |I_2|\,\le \,\left|1+\frac{c}{18}-\frac{c^2}{72}\right|.
|I_2|\,\le \,\left|\upsilon \right|.
\end{equation}
We shall consider cases where
\[\upsilon= 1+\frac{c}{18}-\frac{c^2}{72}\,\geq \,0.\]
From \cite{PS} (see $D_6$ in Lemma~\ref{lem-7}) we deduce that (\ref{eqa}) is satisfied if
\[1+\frac{c}{18}-\frac{c^2}{72}\,\geq \,\frac{1}{12}\left[\left(2+\frac{c}{18}\right)^2\,+\,8\right].\]
This condition is equivalent to $0<c \leq 144/55=2.61818$.
Unfortunately $c=3$ does not lie within this range. But nevertheless one may use Lemma~\ref{lem-7}.
  From (\ref{eq5.2}) and (\ref{eqa}), we get
$$ 40 |\gamma_4| \leq \,c\,+\,\frac{2c^2}{9}\,+\,\frac{c^2}{2}|\upsilon|\, \mbox{ for $c\in(0,2.61818)$}.
$$
Since $\upsilon$ is positive,   we therefore obtain
\begin{equation}\label{eq3ex}
 |\gamma_4| \le \frac{1}{40}\left[c+\frac{c^2}{18}\left(13+\frac{c}{2}-\frac{c^2}{8} \right)\right]
 \, \mbox{ for $c\in(0,2.61818)$}.
 \end{equation}
Further discussion in Lemma~\ref{lem-7} reveals that among the cases therein we have to choose the case $D_9$.
For the case $D_9$, the inequality
$$-\frac{2}{3}\,(|\mu|+1)\,\le \upsilon \,\le \frac{2|\mu|\,(|\mu|+1)}{\mu^2+2|\mu|+4}
$$
is true for $c\in(2.71569,3]$. Clearly the left hand side of the inequality is true for $c\in(0,3]$ whereas the right hand side of the inequality holds only for $c\in(2.71569,3]$. In view of these reasonings, from \eqref{eq1ex}, we find that
\begin{equation}\label{eqb}
|I_2|\le \frac{2}{3}\left(3+\frac{c}{18}\right)
\left[\frac{3+\frac{c}{18}}{3\left(4+\frac{2c}{18}-\frac{c^2}{72} \right)} \right]^{\frac{1}{2}}
= \frac{54+c}{27}\left[\frac{4(54+c)}{3(288+8c-c^2)} \right]^{\frac{1}{2}}
%{\color{red}\mbox{ for $c\in(2.71569,3)$ }}.
\end{equation}
for $c\in(2.71569,3]$.   From (\ref{eq5.2}), we obtain
\begin{equation}\label{eq4ex}
 |\gamma_4| \le  \frac{1}{40}\left(\,c\,+\,\frac{2c^2}{9}\,+\,\frac{c^2}{2}|I_2| \right) \mbox{ for $c\in(2.71569,3]$ }.
 \end{equation}
Hence, from (\ref{eq3ex}) and (\ref{eq4ex}), we obtain the desired inequality.

%In the case  $c=3$,
%$I_2$ from \eqref{eq1ex} reduces to
%$$I_2 = c_3\,+\mu c_1c_2\,+\upsilon c_1^3, \quad \mu =13/6~\mbox{ and }~\upsilon = 25/24.
%$$
% A numerical inspection of the discussion in Lemma~\ref{lem-7} reveals that among the cases therein
% we have to choose the case $D_9$ and as view of this,
%For $c=3$, from (\ref{eqb}), we obtain
%\[ |I_2|\le \frac{19}{9}\sqrt{\frac{76}{303}}.
%\]
For example, if we let $c=3$ in \eqref{eq4ex} and \eqref{eqb}, then \eqref{eq4ex} gives
\[|\gamma_4|\le  \frac{1}{40}\left(3\,+2+\,\frac{19}{2}\sqrt{\frac{76}{303}}\right)
=\frac{1}{40}\left(5\,+\,\frac{19}{2}\sqrt{\frac{76}{303}}\right)\approx 0.243945.\]
On the other hand, for the function $f_3$, we have $\gamma_4(f_3)= \frac{31}{128} \approx 0.2421875$, see (\ref{eq-cor2}).

%%%%%%%%%%%%%%%%%%%%%%%% End \gamma_4 bound %%%%%%%%%%%%%%%%%%%%%%%%
% \section{Fifth logarithmic coefficients bound for $\mathcal{F}(c)$}
Now, we calculate the fifth logarithmic coefficient bound for the class $\mathcal{F}(c)$.
From (\ref{eq2.1}), we compare the coefficients of $z^5$ and get
\[ \beta_5+\beta_4\delta_1+\beta_3\delta_2+\beta_2\delta_3+\beta_1\delta_4+\delta_5
= 2\delta_5+2\delta_1\delta_4+2\delta_2\delta_3+5\delta_5=7\delta_5+2\delta_1\delta_4+2\delta_2\delta_3.
 \]
Using the relations (\ref{eq3}), (\ref{eq4}), (\ref{eq5}) %(\ref{eq5**})
and then simplifying, we get
\begin{align*}
 6\delta_5 & = \beta_5+\beta_4\delta_1+\delta_2(\beta_3-\delta_3)+\delta_3(\beta_2-\delta_2)\\
           & = \beta_5+\frac{1}{2}\beta_1\beta_4+ \frac{1}{24}\beta_1^2\beta_3\,
               +\,\frac{5}{12}\beta_2\left(\beta_3\,+\,\frac{1}{10}\beta_1\beta_2\,
               -\,\frac{1}{20}\beta_1^3\right).
\end{align*}
If we take similar steps as above and use (\ref{eq-s1}), we arrive at the estimation
\begin{align}\label{eq6-s1}
 \nonumber |\gamma_5|  & \le  \frac{1}{60}\left[c\,+\,\frac{c^2}{2}\,+\,\frac{c^3}{24}\,
                +\,\frac{5c^2}{12}\left|c_3\,+\,\left(2+\frac{c}{10}\right)c_1c_2\,
                +\,\left(1+\frac{c}{10}-\frac{c^2}{20}\right)c_1^3\right|\right]\\
             & = \frac{1}{60}\left[c\,+\,\frac{c^2}{2}\,+\,\frac{c^3}{24}\,
                +\,\frac{5c^2}{12}|I_3| \right],
\end{align}
where
\begin{equation}\label{eq2ex}
 I_3:= c_3\,+\mu c_1c_2\,+\upsilon c_1^3,
\quad \mu=2+ \frac{c}{10} ~\mbox{ and }~\upsilon = 1+\frac{c}{10}-\frac{c^2}{20}.
\end{equation}
Our aim is to find an upper bound for $|I_3|$.
Since
\[\upsilon \geq \frac{1}{12}\left[\left(2+\frac{c}{10}\right)^2\,+\,8\right]
% 1+\frac{c}{10}-\frac{c^2}{20}\geq \frac{1}{12}\left(\left(2+\frac{c}{10}\right)^2\,+\,8\right)
\]
for $0<c\leq 80/61 = 1.31148...$, we have proved the sharp estimate for these values of the parameter $c$.
Lemma~\ref{lem-7} gives $|I_3|\le |\upsilon|,\,c\in(0,1.31148)$, for the case $D_6$. Since $\upsilon$ is positive, we therefore
have
\begin{align}\label{eq5ex}
 \nonumber |\gamma_5|  & \le  \frac{1}{60}\left[c\,+\,\frac{c^2}{2}\,+\,\frac{c^3}{24}\,
   +\,\frac{5c^2}{12}\upsilon \right]\\
& = \frac{1}{60}\left[c\,+\,\frac{c^2}{12}\left(11+c- \frac{c^2}{4}\right)\right] \mbox{ for } c\in(1.31148,3].
\end{align}
We see that this corresponds again to the case $D_9$ in Lemma~\ref{lem-7}.
The case $D_9$ holds for $c\in(1.35541,3]$. In view of this, from \eqref{eq2ex}, we get
\begin{equation}\label{eqc}
|I_3|\le \frac{2}{3}\left(3+\frac{c}{10}\right)
\left[\frac{3+\frac{c}{10}}{3\left(4+\frac{2c}{10}-\frac{c^2}{20} \right)} \right]^{\frac{1}{2}}
= \frac{30+c}{15}\left[\frac{2(30+c)}{3(80+4c-c^2)} \right]^{\frac{1}{2}}
\end{equation}
for $c\in(1.35541,3]$. From (\ref{eq6-s1}), we find that
\begin{equation}\label{eq6ex}
  |\gamma_5| \le  \frac{1}{60}\left[c\,+\,\frac{c^2}{2}\,+\,\frac{c^3}{24}\,
                +\,\frac{5c^2}{12}|I_3| \right] \mbox{ for $c\in(1.35541,3]$.}
\end{equation}
Hence, from (\ref{eq5ex}) and (\ref{eq6ex}), we obtain the desired inequality.
For instance, if we take $c=3$, then
%from \eqref{eq2ex}, we find that $\mu =23/10$ and $\upsilon =17/20$.
% We see that this corresponds again to the case $D_9$ in Lemma~\ref{lem-7} and from this
from (\ref{eqc}) we get
\[ |I_3| \le \frac{11}{5}\sqrt{\frac{22}{83}}.
 \]
which by (\ref{eq6-s1}) gives
\[ |\gamma_5|\le  \frac{1}{60}\left(3+\,\frac{9}{2}\,+\,\frac{27}{24}\,
+\,\frac{33}{4}\sqrt{\frac{22}{83}}\right)=\frac{1}{60}\left(\frac{69}{8}\,
+\,\frac{33}{4}\sqrt{\frac{22}{83}}\right)\approx 0.2145050. \]
On the other hand, in this case, $\gamma_5(f_3)= \frac{63}{320}\approx 0.196875$, see (\ref{eq-cor2}).
%%%%%%%%%%%%%%%%%%%%%%%% End \gamma_5 bound %%%%%%%%%%%%%%%%%%%%%%%%
\hfill$\Box$
%%%%%%%%%%%%%%%%%%%%%%%%%%%%%
% \bigskip
% \noindent
\begin{rem}\label{rem-}
% For $\alpha\,=\,1\,-\,c/2\in [0,1)$ the family as $\mathcal{F}(\alpha)$ are known as the families of
% convex functions of order $\alpha$. Hence,
The above estimates contain sharp estimates for the logarithmic
coefficients of convex functions of order $\alpha$.
Especially, we have the sharp estimates for $|\gamma_n|, n=1,2,3,4,$
for convex functions, i.e. for the case $\alpha = 0$, and sharp estimates for
$|\gamma_n|, n=1,2,3,4,5,$ for convex functions  of order $\alpha =1/2$.
\end{rem}
%%%%%%%%%% End- theorem of F(c)  %%%%%%%%%%%%%%%%%%%%%%%%%%%%%%%%%
%%%%%%%%%%%%%%%%%%%%%%%%%%%%%%%%
%%%%%%%%%%%%%%%% Proof of Theorem 6 of G(c)%%%%%%%%%%%%
\subsection{Proof of Theorem~\ref{thm6}}
 Choosing similar abbreviations as in Theorem~\ref{thm5},
$f\in \mathcal{G}(c)$ if and only if
$$ \,1\,+\,\frac{zf''(z)}{f'(z)}\,\prec \,1\,-\,\frac{cz}{1-z}.$$
Hence, we get
\[
\beta_1\,=\,-cc_1,\quad \beta_2\,=\,-c\left(c_2\,+\,c_1^2\right),~\mbox{ and }~
\beta_3\,=\,-c\left(c_3\,+\,2c_1c_2\,+\,c_1^3\right).\]
This results in the relation
\[ \gamma_1\,=\,\frac{\beta_1}{4}\,=\,\frac{-cc_1}{4}
\]
so that $|\gamma_1|\leq c/4,$ where equality is attained for $f'(z)\,=\,(1-z)^c.$
Further,
\[
\gamma_2\,=\,\frac{-c}{12}\left[c_2\,+\,\left(1-\frac{c}{4}\right)c_1^2\right].\]
To get an estimate for $|\gamma_2|$, we use the triangle inequality together with the well-known inequality $|c_2|\leq 1\,-\,|c_1|^2$
for $\phi(z)\,=\,\sum_{k=1}^{\infty}c_kz^k$ in $\mathcal{B}$. This implies
$$|\gamma_2|\,\leq\,\frac{c}{12}\left(1\,-\,\frac{c}{4}|c_1|^2\right)\,\leq\,\frac{c}{12}.$$
Equality is attained here for the function $f'(z)\,=\,(1-z^2)^{c/2}.$
In the case $n=3$, we consider
$$\gamma_3\,=\,\frac{-c}{24}\left[c_3\,+\,\left(2-\frac{c}{2}\right)c_2c_1\,+\,\left(1-\frac{c}{2}\right)
c_1^3\right].$$
% Now we consult again \cite{PS} and we see that the second case is valid which assures that in this case we get the inequality
We see that this corresponds to the case $D_2$ in Lemma \ref{lem-7} and from this we get
$|\gamma_3|\,\leq \,\frac{c}{24}.$
Here, equality is attained for $f'(z)\,=\,(1-z^3)^{c/3}.$
%%%%%%%% End of Theorem 6 of G(c) %%%%%%%%%%%%

\subsection*{Acknowledgements.}
%The authors would like to thank the referee for his/ her comments and useful
%suggestions.
This work was completed while the second author was at IIT Madras for a short period during July-August, 2018.
The work of the first author is supported by Mathematical Research Impact Centric Support of  Department of Science and Technology (DST),
India~(MTR/2017/000367). The second author thanks Science and Engineering Research Board, DST, India,
for its support by  SERB National Post-Doctoral Fellowship
(Grant No. PDF/2016/001274).

%\subsection*{Conflict of Interests}
%The authors declare that there is no conflict of interests regarding the publication of this paper.

 %%%%%%%%%%%%%%%%% References %%%%%%%%%%%%%%%%%%%%%

\end{document}